# New exact solutions of differential equations derived by fractional calculus

F. S. Felber
*Starmark, Inc., P. O. Box 270710, San Diego, CA 92198, USA*

Fractional calculus generalizes the derivative and antiderivative operations $d^n/dz^n$ of differential and integral calculus from integer orders $n$ to the entire complex plane. Methods are presented for using this generalized calculus with Laplace transforms of complex-order derivatives to solve analytically many differential equations in physics, facilitate numerical computations, and generate new infinite-series representations of functions. As examples, new exact analytic solutions of differential equations, including new generalized Bessel equations with complex-power-law variable coefficients, are derived.

This paper presents methods for finding new classes of exact solutions of ordinary differential equations (ODEs). Like the Laplace transform, which is an integral feature, these new methods find exact analytic solutions for whole classes of important ODEs in physics, engineering, and other fields, for which no other means of exact solution are known. To demonstrate the power of these new methods, this paper presents several examples of exact solutions of ODEs with physical significance, some of which are not among the 6200 known exact solutions of ODEs in the most comprehensive and recent handbook [1]. Then the methods are validated by demonstrating in each case that the Runge-Kutta numerical solutions precisely match the new exact analytic solutions. The examples include: (1) a nonhomogeneous Abel equation describing "quasi-exponential" decay with a variable decay rate and a variable source rate; (2) a generalized Bessel equation with variable power-law coefficients describing, for example, radial acoustical eigenmodes of a radially nonuniform circular plate or membrane; and (3) a forced harmonic oscillator equation with a combination of growing and decaying complex-power-law force terms.

These examples barely scratch the surface of the classes of ODEs that can now be solved exactly by these new methods. A big advantage of exact analytic solutions, of course, is that by scaling the solutions as necessary, solutions $F(t)$ can be found immediately for any $t$ without having to integrate numerically from $t = 0$. This advantage can greatly add to accuracy and facilitate certain numerical computations in physics.

The methods presented in this paper for finding new exact solutions of ODEs use fractional calculus, an approach to generalizing derivative and antiderivative operations that is particularly useful for mathematical and computational applications. Generalizing derivatives from integer order to complex order can greatly expand the utility of derivatives and antiderivatives and unleash useful analytic and computational power. This paper provides just a few examples of how this expanded utility can be realized in deriving exact solutions of differential equations in physics, which do not appear in recent compendiums and handbooks [1–3], and in generating new infinite-series representations of functions, which do not appear in tables and handbooks [4–7].

Much attention has been given to the analysis and solution of ODEs and partial differential equations (PDEs) that contain fractional derivatives, which are called fractional differential equations or fractional-order equations. See Podlubny [8], for example, for various methods of solving such equations. Much less attention has been given to the complementary problem of solving integer-order ODEs in terms of fractional derivatives, which are well-defined functions in their own right. That is the objective of this paper, to develop methods for finding new exact analytic solutions of ODEs in terms of fractional derivatives.

An early paper [9] tabulated fractional derivatives of nine elementary functions and the representations in terms of fractional derivatives of seven special functions, including hypergeometric and Bessel functions. Ten years later, Nishimoto [10] solved some simple linear second-order differential equations of the Fuchs type in terms of fractional derivatives. More recently, Kiryakova [11] generalized operators of fractional integration and differentiation and represented generalized functions in terms of fractional derivatives and integrals.

Some overviews of fractional calculus are found in [12–15]. In an historical account, Srivastava and Saxena [16] cites applications of fractional calculus in fluid flow, rheology, diffusive transport, electrical networks, probability, viscoelasticity, electrochemistry of corrosion, and chemical physics. For example, the molecular theory of the viscoelastic properties of polymer solids can be written in terms of fractional calculus [17]. This reference [17] cites numerous others that apply fractional calculus to problems in viscoelasticity. Another historical account of fractional calculus [18] is in the collection of papers on applications in physics edited by Hilfer [19]. Fractional calculus is used in a linearized "Burgers" equation to simulate the frequency-dependent absorption of nonlinear waves in a fluid [20].

The historical account of Ref. [16] also cites applications of fractional calculus in a number of areas of mathematical analysis, including ODEs and PDEs. Seredynska and Hanyga [21] cite numerous mathematical analyses of anomalous diffusion, nonexponential relaxation, viscoelasticity, poroelasticity, and damping in mechanical systems that use fractional derivatives in ODEs and PDEs. In this reference [21], for example, the nonlinear pendulum equation with a fractional-derivative damping term is analyzed.

The particular approach and notation described here to generalizing derivatives to complex order was chosen to be especially useful for mathematical and computational applications. To be most useful, a generalization of derivatives to complex order should give the same results at integer order as ordinary derivatives. The complex-order derivatives defined here do maintain this correspondence. The complex-order

---
*E-mail address:* felber@san.rr.com

derivatives defined here have other useful features in common with integer-order derivatives as well, such as being linear operators. As defined here, complex-order derivative operations on the sine and cosine functions approximate linear phase shifts over certain domains. And as defined here, complex-order derivatives lend themselves to particularly simple and easily applied forms of the Laplace transform of the derivative operator, which is a valuable property for solving differential equations.

Unless specified otherwise, the term "derivative" and the notation $Dz^{(p)} f(z)$ will now be taken to mean the $p$th-order derivative with respect to $z$ of the function $f(z)$. Here, $z$ is a complex variable and $p$ is a complex number. And unless specified otherwise, the term "derivative" will now be taken to include the derivative operator whenever $r = \text{Re}(p)$ is positive, the antiderivative operator whenever $r$ is negative, and the identity operator whenever $p$ is zero.

With a view towards usefulness in applications, the derivative operator is defined to have the following properties:

Scaling:
$$Dz^{(p)} f(cz) = c^p D(cz)^{(p)} f(cz) , \qquad (1)$$

Linear operator:
$$Dz^{(p)}[af(z) + bg(z)] = a Dz^{(p)} f(z) + b Dz^{(p)} g(z) , \qquad (2)$$

Correspondence (for integers $n$):
$$Dz^{(n)} f(z) = d^n f(z)/dz^n , \qquad (3)$$

Inverse and identity:
$$Dz^{(p)} Dz^{(-p)} f(z) = Dz^{(0)} f(z) + G(z) = f(z) + G(z) , \qquad (4)$$

where $a$, $b$, and $c$ are any complex numbers, $f(z)$ and $g(z)$ are any functions of $z$, and $G(z)$ is a gauge function, discussed below.

A derivative operator is not uniquely defined by the properties in Eqs. (1) – (4). With a further view towards usefulness in solving ODEs, the derivative operator is defined to have this additional property:

Derivative operator definition:
$$Dz^{(p)}(z+c)^q = \frac{\Gamma(1+q)}{\Gamma(1+q-p)}(z+c)^{q-p} , \qquad (5)$$

where $q$ is a complex exponent, and $\Gamma$ is the gamma (factorial) function. $\Gamma(z)$ is single-valued and analytic over the entire complex plane, except at the points $z = -m$ ($m = 0, 1, 2, ...$), where it has simple poles with residue $(-1)^m / m!$ [5]. The derivative $Dz^{(p)}(z+c)^q$ is therefore zero at orders $p = 1 + q + m$, consistent with the correspondence property, Eq. (3).

From Eqs. (2) and (5), the derivative of a constant is given by:

Derivative of constant:
$$Dz^{(p)} c = Dz^{(p)}(cz^0) = cz^{-p}/\Gamma(1-p) , \qquad (6)$$

The derivative of a non-zero constant is zero only at orders $p = 1, 2, 3, ...$ .

As defined in Eq. (5), the derivative operator explicitly satisfies the conditions in Eqs. (1) through (4). Figure 1 plots several real-order derivatives and antiderivatives of the real linear function $x$ with respect to $x$, according to the definition in Eq. (5).

If $A(z)$ is a function whose derivative, $Dz^{(-p)} A(z)$, equals $f(z)$, then $A(z) = Dz^{(p)} f(z)$ is the antiderivative of order $p$ of $f(z)$. The real part of the order of an antiderivative, $\text{Re}(p) = r$, is negative. An arbitrary gauge function $G(z)$, a kind of generalized integration constant, can be added to the antiderivative $A(z) = Dz^{(p)} f(z)$, as long as it satisfies $Dz^{(-p)} G(z) = 0$. From Eq. (5), the most general gauge function resulting from the antiderivative operation of order $p$, $Dz^{(p)}$, with $\text{Re}(p) < 0$, is given by
$$G(z) = z^{-p}(c_1 z^{-1} + c_2 z^{-2} + ... + c_k z^{-k}) , \qquad (7)$$
$$\text{for } (r \le -1, \ 1 \le k \le -r) ,$$
where $k$ is the greatest integer less than or equal to $r$. If $r > -1$, then $G(z)$ is zero.

The Laplace transform of the function $(t-t_0)^q u(t-t_0)$ of the real variable $t$ is [5,22]
$$\mathcal{L}\{(t-t_0)^q u(t-t_0)\} = \exp(-t_0 s)\Gamma(1+q) s^{-(1+q)} , \qquad (8)$$
where $t_0 \ge 0$ is a real constant, $s$ is a complex variable, and the unit step function $u(t-t_0)$ equals 1 for $t > t_0$, ½ for $t = t_0$, and 0 for $t < t_0$. In the following, the unit step function will be omitted, but will be understood to apply to Laplace transforms of translated functions. For the Laplace transform in Eq. (8) to exist, the complex exponent $q$ must satisfy $\text{Re}(q) > -1$.

From Eqs. (5) and (8), the Laplace transform of the derivative of $(t-t_0)^q$ is
$$\mathcal{L}\{Dt^{(p)}(t-t_0)^q\} = s^p \mathcal{L}\{(t-t_0)^q\}$$
$$= \exp(-t_0 s)\Gamma(1+q) s^{-(1+q-p)} . \qquad (9)$$
For this Laplace transform to exist, the complex order $p$ must satisfy $\text{Re}(p) < \text{Re}(q+1)$.

An example of a differential equation with variable complex coefficients that can be solved exactly in terms of complex-order derivatives is
$$dF(t)/dt + ct^a F(t) = t^b , \qquad (10)$$
where $a$, $b$, and $c$ are complex constants. This equation, a nonhomogeneous Abel equation of the first kind, could represent a quantity $F(t)$ that decays in time $t$ at a variable rate $ct^a$, and that is driven by a variable source term $t^b$.

If $a = b$, the solution of Eq. (10) is trivial. If $a = 0$, then Eq. (10) can be transformed to an Emden-Fowler equation that is exactly solvable for arbitrary $b$ [1,2]. In the general case, in which $a$ and $b$ are any complex constants, the solution of Eq. (10) proceeds as follows.

After a change to the new independent variable $w = ct^{a+1}/(a+1)$, which is the integral of the decay rate, the Laplace transform of Eq. (10) becomes
$$f(s) = \frac{F(0)}{s+1} + c^{p-1}(a+1)^{-p}\Gamma(1-p)\frac{s^p}{s(s+1)} , \qquad (11)$$
where $p = (a-b)/(a+1)$, and $f(s) = \mathcal{L}\{F(w)\}$ is the image function of $F(w)$. From Eq. (9), the inverse Laplace transform of Eq. (11) gives the exact analytic solution of Eq. (10) as

$$F(w) = F(0)\exp(-w) + c^{p-1}(a+1)^{-p}\Gamma(1-p)Dw^{(p)}[1-\exp(-w)] \quad . \tag{12}$$

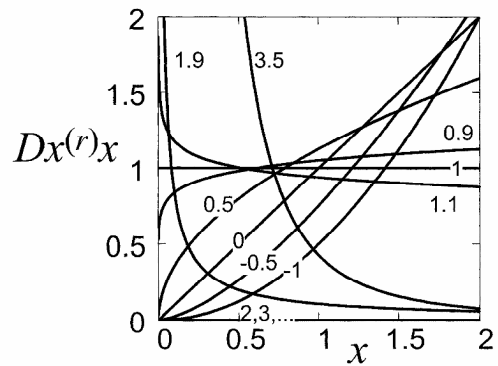

FIG. 1  Derivatives of $x$, $Dx^{(r)}x$, vs. $x$ from Eq. (5) for the orders $r$ indicated.

A typical solution is plotted in Fig. 2.

An example of a second-order differential equation with variable complex coefficients that can be solved exactly in terms of complex-order derivatives is

$$\frac{d^2F(t)}{dt^2} + \frac{a}{t}\frac{dF(t)}{dt} + b^2 t^{2c-2} F(t) = 0 \quad , \tag{13}$$

where $a$, $b$, and $c$ are complex constants. If the variable $t$ were a radius, the solution of Eq. (13) could be, for example, an azimuthally symmetric radial acoustic eigenmode of a circular membrane or thin plate with certain radial nonuniformities of the membrane or plate properties [23]. If $a = b = c = 1$, then Eq. (13) becomes the Bessel equation of zero order, and the solution becomes the zero-order Bessel function of the first kind, $J_0(t)$. If $c = 0$, then Eq. (13) becomes the Euler equation [1,2]. If $2c$ equals a nonzero integer $n$, then the exact solution of Eq. (13) is known in terms of Bessel functions of the first kind of order $(1-a)/n$ [1]. In the general case, the solution of Eq. (13) proceeds as follows.

After a change to the new independent variable $w = (bt^c/2c)^2$, the Laplace transform of Eq. (13) becomes

$$f(s) = F(0)\sum_{m=0}^{\infty}\frac{\Gamma(p-m)}{\Gamma(p)s^{m+1}} \quad , \tag{14}$$

where $p = (1-a)/2c$, and $f(s)$ is the image function of $F(w)$. Using Eq. (9) and $\mathcal{L}\{J_0(2\sqrt{w})\} = s^{-1}\exp(-1/s)$, the inverse Laplace transform of Eq. (14) gives the exact analytic solution of Eq. (13) as

$$F(w) = F(0)\Gamma(1-p)w^p Dw^{(p)} J_0(2\sqrt{w}) \quad . \tag{15}$$

A typical solution is plotted in Fig. 2. Non-integer exponents were used in Fig. 2 to demonstrate that the analytic solutions were found without conventional recursion methods.

In the examples above, after suitable changes of variables, Laplace transforms were used to derive solutions of differential equations in terms of complex-order derivatives. The inverse process can be used as well to generate whole classes of differential equations for which exact solutions in terms of complex-order derivatives exist. The recipe for this inverse process of generating equations from solutions is as follows: One starts with a "solution" involving a complex-order derivative of a well-defined function. Then the inverse Laplace transform of the image function of the "solution" is expressed as a differential equation. The differential equation so derived can be transformed by a change of the independent variable to derive classes of equations for which the "solution" is exact.

As a simple example of this inverse process for generating exactly solvable differential equations, consider the "solution,"

$$F(w) = Dw^{(p)}[a_n\cos(nw) + b_n\sin(nw)] \quad , \tag{16}$$

where $p$ is the complex order of the derivative, and $a_n$, $b_n$, and $n$ are complex constants. If $n$ is a positive integer, then $a_n\cos(nw) + b_n\sin(nw)$ can represent one component of a Fourier series expansion of a periodic function. From Eq. (9), the image function of this "solution" is $f(s) = (a_n s + b_n n)s^p/(s^2 + n^2)$. The inverse Laplace transform generates the forced harmonic oscillator equation,

$$\frac{d^2F(w)}{dw^2} + n^2 F(w) = \frac{-a_n(1+p) + b_n nw}{\Gamma(-p)w^{2+p}} \quad , \tag{17}$$

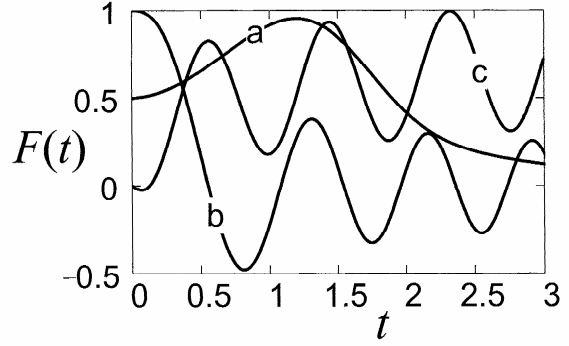

FIG. 2. Exact analytic solutions $F(t)$ of: a) $dF/dt + 0.7t^{3.1}F = t^{0.8}$, $F(0) = 0.5$ from Eq. (12); b) $d^2F/dt^2 + 0.7t^{-1}dF/dt + 6.3^2 t^{0.6}F = 0$, $F(0) = 1$ from Eq. (15); and c) $d^2F/dt^2 + 7.1^2 F = -0.34t^{-0.79} + 27t^{0.21}$ from Eq. (16).

for which Eq. (16) is the exact solution. A typical solution (for $p = -1.21$, $n = 7.1$, $a_n = -1.5$, $b_n = 3.5$) is shown in Fig. 2.

The numerical solutions by the Runge-Kutta technique [24] of Eqs. (10), (13), and (17) match the exact analytic solutions in Fig. 2 precisely. By scaling exact analytic solutions in accordance with Eq. (1), as necessary, $F(t)$ can be found immediately for any $t$ without having to integrate numerically from $t = 0$.

A transformation of the independent variable, $w = vt^c$, in Eq. (17), in which $v$ and $c$ are complex constants, generates the new differential equation,

$$\frac{d^2F(t)}{dt^2} + \frac{(1-c)}{t}\frac{dF(t)}{dt} + (ncv)^2 t^{2c-2} F(t) = \frac{c^2[-a_n(1+p) + b_n nvt^c]}{\Gamma(-p)v^p t^{2+cp}} \quad , \tag{18}$$

for which Eq. (16), after transformation, is also an exact solution. In the same manner, other transformations of the variable in Eq. (17) can generate other classes of ODEs for which Eq. (16), after transformation, is an exact solution.

From Eqs. (1), (2), (5), (9), and the Maclaurin series expansions of the circular functions, the Laplace transforms of the derivatives of the circular functions are

$$\mathcal{L}\{Dt^{(p)}\sin[\omega t - (2\pi j - \theta)]\} = \exp[-(2\pi j - \theta)s]\omega s^p/(s^2 + \omega^2) \quad , \tag{19}$$

$$\mathcal{L}\{Dt^{(p)}\cos[\omega t - (2\pi j - \theta)]\} = \exp[-(2\pi j - \theta)s]s^{p+1}/(s^2 + \omega^2) \quad , \tag{20}$$

where $\omega$ and $\theta$ are real constants, and $j$ is any integer satisfying $2\pi j - \theta \geq 0$.

The inverse Laplace transforms of the image functions in Eqs. (19) and (20) are found from the Bromwich integrals and the calculus of residues [22] to be

$$Dt^{(p)}\sin[\omega t - (2\pi j - \theta)] = \omega^p \sin[\omega t - (2\pi j - \theta) + p\pi/2] \quad , \tag{21}$$

$$Dt^{(p)}\cos[\omega t - (2\pi j - \theta)] = \omega^p \cos[\omega t - (2\pi j - \theta) + p\pi/2] \quad . \tag{22}$$

From Eqs. (5), (21), and (22), therefore, we find that derivatives of the circular functions and their Maclaurin series expansions approach phase-shifted circular functions for sufficiently large arguments $x$, as

$$Dx^{(p)} \sin x = \sum_{m=0}^{\infty} \frac{(-1)^m x^{2m+1-p}}{\Gamma(2m+2-p)} \sim \sin(x + p\pi/2) \; , \qquad (23)$$

$$Dx^{(p)} \cos x = \sum_{m=0}^{\infty} \frac{(-1)^m x^{2m-p}}{\Gamma(2m+1-p)} \sim \cos(x + p\pi/2) \; . \qquad (24)$$

The phase-shifted circular functions on the right-hand sides of Eqs. (23) and (24) approach the derivatives on the left in the asymptotic limit $x \to \infty$. The inexactness of the relationships in Eqs. (23) and (24) at small $x$ seems to be related to the Gibbs phenomenon, which occurs with Fourier series and other eigenfunction expansions, particularly at discontinuities [25]. Figure 3 shows the real domains over which the phase-shifted circular functions in Eqs. (23) and (24) are "good approximations," defined in Fig. 3 as the domains in which the right-hand and left-hand sides of Eqs. (23) and (24) differ from each other by no more than 0.01 for the first 40 terms in the series. Since the complex-order derivatives of the circular functions are equivalent to phase shifts proportional to the order, at least over reasonably useful domains, the choice of derivative operator definition in Eq. (5) has enhanced utility for physics and computational applications.

Many new series expansions can also be generated from complex-order derivatives. For example, Eqs. (16) and (17) give a rapidly converging alternate-series expansion of the constant function 1,

$$1 = \Gamma^n(c) \sum_{m=0}^{\infty} (-1)^m \left[ \frac{w^{2m}}{\Gamma^n(2m+c)} + \frac{w^{2m+2}}{\Gamma^n(2m+2+c)} \right], \qquad (25)$$

where $w$ is a complex variable, $c$ is any complex constant, and $n$ is any positive integer. The convergence to 1 of Eq. (25) can be improved without limit by choosing $n$ arbitrarily high.

In summary, an approach was presented for generalizing derivatives and antiderivatives from integer orders to complex orders, preserving some, but not all, important properties of integer-order calculus. The definition of a complex-order derivative operator was chosen to have useful computational applications in physics. The choice of operator definition makes Laplace transforms particularly useful for solving differential equations. An inverse process was also presented for generating classes of differential equations that are exactly solvable with complex-order derivatives. These new solution methods will substantially expand the set of exactly solvable differential equations. As examples, this paper presented exact solutions in terms of complex-order derivatives of a nonhomogeneous decay equation with a variable decay rate, a generalized Bessel equation, and a forced harmonic-oscillator equation. These few examples of new exact analytic solutions of differential equations and new series expansions of functions found with complex-order derivatives offer only a small glimpse of the potential utility of fractional calculus for mathematical applications and computations.

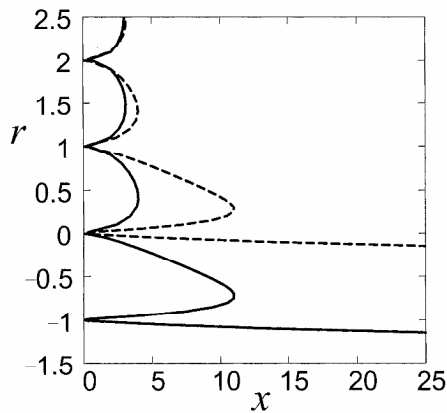

FIG. 3. Left boundaries of semi-infinite domains $x$ in which $\cos(x + r\pi/2)$ approximates (within 0.01) $Dx^{(r)} \cos x$ (solid curves) and $\sin(x + r\pi/2)$ approximates $Dx^{(r)} \sin x$ (dashed curves) for real order $r$.